\documentclass[12pt]{amsart}
\usepackage{amssymb,amsmath}
\numberwithin{equation}{section}

\def\om{\omega}
\def\bom{\bar\omega}
\def\H{\mathcal H}
\def\sumj{\underset{j=1}{\overset {n-1}\sum}}
\def\dib{\bar\partial}
\def\di{\partial}
\def\simleq{\underset\sim<}

\def\T{\text}
\def\1#1{\overline{#1}}
\def\2#1{\widetilde{#1}}
\def\3#1{\widehat{#1}}
\def\4#1{\mathbb{#1}}
\def\5#1{\frak{#1}}
\def\6#1{{\mathcal{#1}}}
\def\C{{\4C}}
\def\R{{\4R}}

\def\sumK{\underset{|K|=k-1}{{\sum}'}}
\def\sumJ{\underset{|J|=k}{{\sum}'}}
\def\sumij{\underset {ij=1,\dots,n-1}{{\sum}}}

\def\sumjq{\underset {j\leq q}\sum}

\def\B{\mathcal B}

\def\phi{\varphi}

\newtheorem{Thm}{Theorem}[section]
\newtheorem{Cor}[Thm]{Corollary}
\newtheorem{Pro}[Thm]{Proposition}
\newtheorem{Lem}[Thm]{Lemma}
\theoremstyle{definition}\newtheorem{Def}[Thm]{Definition}
\theoremstyle{remark}
\newtheorem{Rem}[Thm]{Remark}
\newtheorem{Exa}[Thm]{Example}

\def\Label#1{\label{#1}}
\def\bl{\begin{Lem}}
\def\el{\end{Lem}}
\def\bp{\begin{Pro}}
\def\ep{\end{Pro}}
\def\bt{\begin{Thm}}
\def\et{\end{Thm}}
\def\bc{\begin{Cor}}
\def\ec{\end{Cor}}
\def\bd{\begin{Def}}
\def\ed{\end{Def}}
\def\br{\begin{Rem}}
\def\er{\end{Rem}}
\def\be{\begin{Exa}}
\def\ee{\end{Exa}}
\def\bpf{\begin{proof}}
\def\epf{\end{proof}}
\def\ben{\begin{enumerate}}
\def\een{\end{enumerate}}

\def\1alpha{[\frac1\alpha]}
\def\T{\text}
\def\R{{\Bbb R}}

\def\C{{\Bbb C}}

\numberwithin{equation}{section}
\def\T{\text}

\newcommand{\no}[1]{\|{#1}\|}
\newcommand{\NO}[1]{\|{#1}\|^2}

\newtheorem{theorem}{Theorem  }[section]
\newtheorem{definition}[theorem]{Definition }
\newtheorem{lemma}[theorem]{Lemma  }
\newtheorem{proposition}[theorem]{Proposition  }
\newtheorem{corollary}[theorem]{Corollary }
\newtheorem{example}[theorem]{\it Example }

\begin{document}

\title[Compactness of $\Box_b$ on a CR manifold]{Compactness of $\Box_b$ on a CR manifold}         
\author[ T.V.~Khanh, S.~Pinton and G.~Zampieri ]
{Tran Vu Khanh, Stefano Pinton and Giuseppe Zampieri}
\address{Faculty of Mathematics and Computer Science, University of Science - Ho Chi Minh City, Vietnam}
\address{Dipartimento di Matematica, Universit\`a di Padova, via 
Trieste 63, 35121 Padova, Italy}
\email{tvkhanh@math.hcmuns.edu.vn}
\email{pinton@math.unipd.it}
\email{ zampieri@math.unipd.it}     
\maketitle
\begin{abstract}
This note is aimed at simplifying current literature about compactness estimates for the Kohn-Laplacian on CR manifolds. The approach consists in a tangential basic estimate in the formulation given by the first author in \cite{Kh10} which refines former work by Nicoara \cite{N06}. 
It has been proved by Raich \cite{R10} that on a CR manifold of dimension $2n-1$ which is compact pseudoconvex of hypersurface type embedded in $\C^n$ and orientable, the property named ``$(CR-P_q)$" for $1\leq q\leq \frac{n-1}2$, a generalization of the one introduced by Catlin in \cite{C84}, implies compactness estimates for the Kohn-Laplacian $\Box_b$ in degree $k$ for any $k$ satisfying $q\leq k\leq n-1-q$. The same result is stated by Straube in \cite{S10} without the assumption of orientability.
We regain these results by a simplified method and extend the conclusions in two directions. First, the CR manifold is no longer required to be embedded. Second, when $(CR-P_q)$ holds for $q=1$ (and, in case $n=1$, under the additional hypothesis that $\dib_b$ has closed range on functions) we prove compactness also in the critical degrees $k=0$ and $k=n-1$.

\end{abstract}
\tableofcontents 
\section{Introduction and statements}
Let $M$ be a compact pseudoconvex CR manifold of hypersurface type of real dimension $2n-1$ endowed with the Cauchy-Riemann structure $T^{1,0}M$. We choose a basis $L_1,...,L_{n-1}$ of $T^{1,0}M$, the conjugated basis $\bar L_1,...,\bar L_{n-1}$ of $T^{0,1}M$, and a transversal, purely imaginary, vector field $T$. We also take a hermitian metric on the complexified tangent bundle in which we get an orthogonal decomposition $\C TM=T^{1,0}M\oplus T^{0,1}M\oplus \C T$. We denote by $\om_1,...,\om_{n-1},\bom_1,...,\bom_{n-1},\gamma$ the dual basis of $1$-forms. We denote by $\mathcal L_M$ the Levi form defined by $\mathcal L_M(L,\bar L'):=d\gamma(L,\bar L') $ for $L,\,L'\in T^{1,0}M$. The coefficients of the matrix $(c_{ij})$ of $\mathcal L_M$ in the above basis are described through Cartan formula as
$$
c_{ij}=\langle \gamma,[L_i,\bar L_j]\rangle.
$$
We denote by $\6B^k$ the space of $(0,k)$-forms $u$ with $C^\infty$ coefficients; they are expressed, in the local basis, as $u=\sumJ u_J\bom_J$ for $\bom_
J=\bom_{j_1}\wedge...\wedge \bom_{j_k}$. Associated to the Riemaniann metric $\langle\cdot,\cdot\rangle_z,\,\,z\in M$ and to the element of volume $dV$, 
there is a $L^2$-inner product $(u,v)=\int_M\langle u,v\rangle_z dV$. We denote by $(L^2)^k$ the completion of $\6B^k$ under this norm; we also use the 
notation $(H^s)^k$ for the completion under the Sobolev norm $H^s$. Over the spaces $\6B^k$ there is induced by the de-Rham exterior derivative a complex $\dib
_b:\,\6B^k\to \6B^{k+1}$. We denote by $\dib_b^*:\,\6B^k\to\6B^{k-1}$ the adjoint and set $\Box_b=\dib_b\dib^*_b+\dib^*_b\dib_b$. Let $\phi$ be a smooth 
function and denote by $(\phi_{ij})$ the matrix of the Levi form $\6L_\phi=\frac12(\di_b\dib_b-\dib_b\di_b)(\phi)$ in the  basis above. Let $L^2_\phi$ be 
the $L^2$ space weighted by $e^{-\phi}$ and denote by $L_j^\phi=L_j-\phi_j$, for $\phi_j:=L_j(\phi)$, the $L^2_\phi$-adjoint of $-\bar L_j$. The following  is 
the tangential version of the celebrated H\"ormander-Kohn-Morrey basic estimate. We present here the refinement by Khanh \cite{Kh10} of a former statement 
by Nicoara \cite{N06}. Le $z_o\in M$; for a suitable neighborhood $U$ of $z_o$ and a constant $c>0$, we have
\begin{eqnarray} \Label{hkm}
\begin{split}
\no{\bar{\partial}_b u}^2_{\phi}&+\no{\bar{\partial}^*_{b,\phi} u}^2_{\phi}+c\no{u}^2_{\phi}\\
\ge & {\sumK}\sum_{ij}(\phi_{ij}u_{iK},u_{jK})_\phi-{\sumJ}\sum_{j=1}^{q_o}(\phi_{jj}u_J,u_J)_\phi\\
&+{\sumK}\sum_{ij}(c_{ij}Tu_{iK},u_{jK})_\phi-{\sumJ}\sum_{j=1}^{q_o}(c_{jj}Tu_J,u_J)_\phi\\
&+\frac{1}{2}\big(\sum^{q_o}_{j=1}\no{L_j^{\phi} u}^2_{\phi}+\sum^{n-1}_{j=q_o+1}\no{\bar{L}_ju}^2_{\phi} \big),
\end{split}
\end{eqnarray}
for any $u\in \B^{k}_c(U)$ where $q_o$ is any integer  with $0\le q_o\le n-1$. We introduce now a potential-theoretical condition which is a variant 
of the ``$P$-property" by Catlin \cite{C84}. In the present version it has been introduced by Raich \cite{R10} and Straube \cite{S10}.
\bd
\Label{d1.1}
Let $z_o$ be a point of $M$ and $q$ an index in the range $1\leq q\leq n-1$. We say that $M$ satisfies property $(CR-P_q) $ at $z_o$ if there is a family of weights $\{\phi^\epsilon\}$ in a neighborhood $U$ of $z_o$ such that, if $\lambda_1^{\phi^\epsilon}\leq...\leq \lambda^{\phi^\epsilon}_{n-1}$ are the ordered eigenvalues of the Levi form $\6L_{\phi^\epsilon}$, we have
\begin{equation}
\Label{supernova}
\begin{cases}
|\phi^\epsilon(z)|\leq 1,\quad &z\in U
\\
\underset{j=1}{\overset q\sum}\lambda^{\phi^\epsilon}_j(z)\geq \epsilon^{-1},\quad&z\in U\T{ and $\ker\mathcal L_M(z)\neq\{0\}$.}
\end{cases}
\end{equation}
\ed
It is obvious that $(CR-P_q)$ implies $(CR-P_k)$ for any $k\geq q$.
\br
\Label{r1.0}
We restrict our considerations to the unit sphere. Outside a neighborhood $V_\epsilon$ of $\ker d\gamma$, the sum $\sumjq \lambda_j^{\phi^\epsilon}$ can get negative; let $-b_\epsilon$ be a bound from below. Now, if $c_\epsilon$ is a bound from below for $d\gamma$ outside  $V_\epsilon$, by setting $a_\epsilon:=\frac{\epsilon^{-1}+b_\epsilon}{qc_\epsilon}$, we have, 
\begin{equation}
\Label{1.0}
\sumjq\lambda^{\phi^\epsilon}_j+qa_\epsilon c_\epsilon\geq\epsilon^{-1}\quad\T{on the whole $U$}.
\end{equation}
Again, \eqref{1.0} for $q$ implies \eqref{1.0} for any $k\geq q$. 
\er
Related to this notion there is the main result of the paper
\bt
\Label{t1.1}
Let $M$ be a compact pseudoconconvex CR manifold of hypersurface type of dimension $2n-1$. Assume that $(CR-P_q)$ holds for $1\leq q\leq \frac{n-1}2$ over a covering $\{U\}$ of $M$. Then we have compactness estimates
\begin{equation}
\Label{1.2}
\NO{u}\leq \epsilon(\NO{\dib_b u}+\NO{\dib^*_bu})+C_\epsilon\NO{u}_{-1}\quad\T{for any $u\in D^k_{\dib^*_b}\cap D^k_{\dib_b} $ with $q\leq k\leq n-1-q$,}
\end{equation}
\et
where $D^k_{\dib^*_b}$ and $D^k_{\dib_b}$ are the domains of $\dib^*_b$ and $\dib_b$ respectively. 
The proof of this, as well as of the theorem which follows, is given in Section~\ref{s2}.
Let $\H^k=\ker \dib_b\cap\ker\dib_b^*$ be the space of harmonic forms of degree $k$. As a consequence of \eqref{1.2}, we 
have that for $q\leq k\leq n-1-q$, 
the space $\H^k$ is finite-dimensional, $\Box_b$ is invertible over $\H^{k\,\perp}$ (cf. \cite{N06} Lemma 5.3) and its inverse $G_k$ is a compact operator. When $k
=0$ and $k=n-1$ it is no longer true that it is finite-dimensional. However, if $q=1$, we have a result analogous to \eqref{1.2} also in the critical 
degrees $k=0$ and $k=n-1$.
\bt
\Label{t1.2}
Let $M$ be a compact, pseudoconvex CR manifold of hypersurface type of dimension $2n-1$. Assume that property $(CR-P_q)$ holds for $q=1$ over a covering $\{U\}$ of $M$ and, in case $n=2$, make the additional hypothesis that $\dib_b$ has closed range. Then
\begin{equation}
\NO{u}\leq \epsilon(\NO{\dib_b u}+\NO{\dib^*_bu})+C_\epsilon\NO{u}_{-1}\quad \T{for any $u\in\H^{k\,\perp},\,\,k=0$ and $k=n-1$.}
\end{equation}
In particular, $G_k$ is compact for $k=0$ and $k=n-1$.
\et

\section{Proofs}
\Label{s2}
{\it Proof of Theorem~\ref{t1.1}.}\hskip0.3cm
We choose a local patch where a local frame of vector fields is found in which \eqref{hkm} is fulfilled.
The key point is to specify the convenient choices of $q_o$ and $\phi$ in \eqref{hkm}. Let $1=\psi^++\psi^-+\psi^0$ be a conic, smooth partition of the unity in  $\R^{2n-1}$ dual to the space to which $U$  is identified in local coordinates. Let $\T{id}=\Psi^++\Psi^-+\Psi^0$ be the microlocal decomposition of the identity  by the pseudodifferential operators with symbols $\psi^{\overset\pm0}$ and let $\zeta$ be a cut-off function. We decompose a form $u$ as
\begin{equation}
\Label{2.1}
u^{\overset\pm0}=\zeta\Psi^{\overset\pm0}u\quad u\in\B^k_c(U),\,\, \zeta|_{\T{supp}\,u}\equiv1.
\end{equation}
For $u^+$ we choose $q_o=0$ and $\phi=\phi^\epsilon$. 
We also need to go back to Remark~\ref{r1.0}. Now, if $a_\epsilon$ has been chosen so that \eqref{1.0} is fulfilled, we remove $T$ from our scalar products observing that, for large $\xi$, we have $\xi_{2n+1}>a_\epsilon$ over $\T{supp}\,\psi^+$. In the same way as in Lemma 4.12 of \cite{N06}, we conclude that for $k\geq q$ 
\begin{equation*}
\begin{split}
\sumK\sumij&(( c_{ij}T+\phi^\epsilon_{ij})u_{iK}^+,u^+_{jK})_{\phi^\epsilon}\geq \sumK\sumij((a_\epsilon c_{ij}+\phi^\epsilon_{ij})u_{iK}^+,u^+_{jK})_{\phi^\epsilon}\\
& -C\NO{u^+}_{\phi^\epsilon}-C_\epsilon\NO{u^+}_{-1,\phi^\epsilon}
-C_\epsilon\NO{\tilde\Psi^0u^+}_{\phi^\epsilon}\\
&\geq\epsilon^{-1}\NO{u^+}_{\phi^\epsilon}-C_\epsilon\NO{u^+}_{-1,\phi^\epsilon}-C_\epsilon\NO{\tilde\Psi^0u^+}_{\phi^\epsilon},
\end{split}
\end{equation*}
where $\tilde\Psi^0\succ\Psi^0$ in the sense that $\tilde \psi^0|_{supp\,\psi^0}\equiv1$.
Note that there is here an additional term $-C_\epsilon\NO{u^+}_{-1,\phi^\epsilon}$ with respect to \cite{N06}. Reason is that $(c_{ij}\xi_{2n-1}+\phi^\epsilon_{ij})$ can get negative values, even on $\T{supp}\,\psi^+$, when $\xi_{2n-1}<a_\epsilon$. Integration in this compact region, produces the above error term.
 It follows
\begin{equation}
\Label{2.2}
\NO{u^+}_{\phi^\epsilon}\leq \epsilon(\NO{\dib_bu^+}_{\phi^\epsilon}+\NO{\dib_{b,\phi^\epsilon}^*u^+}_{\phi^\epsilon})+C_\epsilon\NO{u^+}_{-1,\phi^\epsilon}+C_\epsilon\NO{\tilde\Psi^0u^+}_{\phi^\epsilon},\quad k=1,...,n-1.
\end{equation}
By taking composition $\chi(\phi^\epsilon)$ for a smooth function $\chi=\chi(t),\,t\in\R^+$ with $\dot\chi>0$ and $\ddot\chi>0$, we get
$$
(\chi(\phi^\epsilon))_{ij}=\dot\chi \phi^\epsilon_{ij}+\ddot\chi|\phi^\epsilon_j|^2\kappa_{ij},
$$
where $\kappa_{ij}$ is the Kronecker symbol.
We also notice that
$$
|\dib^*_{b,\chi(\phi^\epsilon)}u|^2\leq 2|\dib_b^*u|^2+2\dot\chi^2\sumK|\sumj\phi^\epsilon_ju_{jK}|^2.
$$
Remember that $\{\phi^\epsilon\}$ are uniformly bounded by $1$. Thus, if we choose $\chi=e^{(t-1)}$, we have that $\ddot\chi\geq 2\dot\chi^2$ for $t=\frac12\phi^\epsilon$. For this reason, with this modified weight, we can replace the weighted adjoint $\dib^*_{b,\phi^\epsilon}$ by the unweighted $\dib_b^*$ in \eqref{2.2}. By the uniform boundedness of the weights, we can also remove them from the norms and end up with the estimate
\begin{equation}
\Label{2.3}
\epsilon^{-1}\NO{u^+}\leq \NO{\dib^*_bu^+}+\NO{\dib_bu^+}+C_\epsilon\NO{\tilde \Psi^0u},\quad k=q,...,n-1.
\end{equation}
For $u^-$, we choose $q_o=n-1$ and $\phi=-\phi^\epsilon$. Observe that for $|\xi| $ large we have  $-\xi_{2n-1}\geq a_\epsilon$ over $\T{supp}\,\psi^-$ (cf. \cite{N06} Lemma 4.13); thus, we have in the current case, for $k\leq n-1-q$  
\begin{equation*}
\begin{split}
\sumK&\sumij ((c_{ij}T-\phi^\epsilon_{ij})u^-_{iK},u^-_{jK})_{-\phi^\epsilon}-\sumJ\sumj ((c_{jj}T-\phi_{jj}^\epsilon)u^-_J,u^-_J)_{-\phi^\epsilon}
\\&\geq-\sumK\sumij ((a_\epsilon c_{ij}+\phi^\epsilon_{ij})u^-_{iK},u^-_{jK})_{-\phi^\epsilon}
\\
&+\sumJ\sumj ((a_\epsilon c_{jj}+\phi_{jj}^\epsilon)u^-_J,u^-_J)_{-\phi^\epsilon} 
-C\NO{u^-}_{\phi^\epsilon}-C_\epsilon\NO{u^-}_{-1,\phi^\epsilon}-C_\epsilon\NO{\tilde\Psi^0u^-}_{\phi^\epsilon}
\\
&\geq \epsilon^{-1}\NO{u^-}_{\phi^\epsilon}-C\NO{u^-}_{\phi^\epsilon}-C_\epsilon\NO{u^-}_{-1,\phi^\epsilon}-C_\epsilon\NO{\tilde\Psi^0u^-}_{\phi^\epsilon}.
\end{split}
\end{equation*}
Thus, we get the analogous of \eqref{2.2} for $u^+$ replaced by $u^-$ and, removing again the weight from the adjoint $\dib^*_{b,\phi^\epsilon}$ and from the norms, we conclude
\begin{equation}
\Label{2.4}
\NO{u^-}\leq \epsilon(\NO{\dib_b u^-}+\NO{\dib^*_b u^-})+C_\epsilon\NO{u^-}_{-1,\phi^\epsilon}+C_\epsilon\NO{\tilde \Psi^0u},\quad k=0,...,n-1-q.
\end{equation}
In addition to \eqref{2.3} and \eqref{2.4}, we have elliptic estimates for $u^0$
\begin{equation}
\Label{2.5}
\NO{u^0}_1\simleq \NO{\dib u^0}+\NO{\dib^*_bu^0}+\NO{u}_{-1}.
\end{equation}
We put together \eqref{2.3}, \eqref{2.4} and \eqref{2.5} and notice that 
\begin{equation}
\Label{2.6}
\begin{split}
\NO{\dib_b(\zeta\Psi^{\overset\pm0}u)}&\leq \NO{\zeta\Psi^{\overset\pm0}\dib_bu}+\NO{[\dib_b,\zeta\Psi^{\overset\pm0}]u}
\\
&\leq \NO{\Psi^{\overset\pm0}\dib_bu}+\NO{\tilde\zeta\Psi^{\overset\pm0}u}+\NO{\tilde\Psi^0u},
\end{split}
\end{equation}
for $\tilde \zeta\succ\zeta$ and $\tilde\Psi^0\succ\Psi^0$. The similar estimate holds for $\dib_b$ replaced by $\dib_b^*$. Since $\zeta|_{\T{supp}\,u}\equiv1$, then
\begin{equation*}
\begin{split}
\NO{u}&\leq\underset{+,-,0}\sum\NO{\zeta\Psi^{\overset\pm0}u}+Op^{-\infty}(u)
\\
&\leq\epsilon\underset{+,-,0}\sum(\NO{(\dib_bu)^{\overset\pm0}}+\NO{(\dib_b^*u)^{\overset\pm0}})+C_\epsilon\NO{u}_{-1},
\end{split}
\end{equation*}
and therefore
\begin{equation}
\Label{2.7}
\NO{u}\leq \epsilon(\NO{\dib_bu}+\NO{\dib_b^*u})+C_\epsilon\NO{u}_{-1},\quad q\leq k\leq n-1-q.
\end{equation}
We pass now to consider $u$ globally defined on the whole $M$ instead of a local patch $U$. We cover $M$ by $\{U_\nu\}$ so that in each patch there is a basis of forms in which the basic estimate holds. In the identification of $U_\nu$ to $\R^{2n-1}$, we suppose that the microlocal decomposition which yields \eqref{2.7} is well defined. We then apply \eqref{2.7} to a decomposition $u=\sum_\nu \zeta_\nu u$ for a partition of the unity $\sum_\nu\zeta_\nu=1$ on $M$, observe that $[\dib_b,\zeta_\nu]$ and  $[\dib_b^*,\zeta_\nu]$ are $0$-order operators and thus give rise to an error term which is controlled by $\epsilon^{-1}\NO{u}$ and get \eqref{2.7} for any $u\in \B^k$. Finally,  by the density of smooth forms $\B^k$ into Sobolev forms  $(H^1)^k$, \eqref{2.7} holds in fact for any $u\in D^k_{\dib^*_b}\cap D^k_{\dib_b}$. The proof is complete.

\hskip15cm$\Box$

\vskip0.3cm
{\it Proof of Theorem~\ref{t1.2}.}\hskip0.3cm
When $q=1$, we observe that we have the estimate for $u^-$ in degree $k=0$ and for $u^+$ in degree $k=n-1$ (cf. \cite{KN06} Lemma 3.3). We prove how to get the estimate for $u^+$ in degree $k=0$ (the one for $u^-$ when $k=n-1$ being similar). Now, if $n>2$, we have, as a consequence of \eqref{2.7}, that $\dib_b^*$ has closed range on $1$-forms. In particular
\begin{equation*}
\begin{split}
\H^{0\,\perp}&=(\ker\dib_b)^\perp
\\&=\T{range}\,\dib_b^*.
\end{split}
\end{equation*}
Thus, if $u\in\H^{0\,\perp}$, then $u=\dib^*_bv$ for some $v\in (L^2)^1$. Moreover, we can choose the solution $v$ belonging to $(\ker \dib_b^*)^\perp\subset\overline{\T{range}\,\dib_b}\subset \ker\dib_b$. Also, since $v\in(\ker \dib_b^*)^\perp\subset\mathcal H^{0\,\perp}$, then (cf. \cite{N06} Lemma 5.3)
\begin{equation*}
\begin{split}
\NO{v}_0&\simleq\NO{\dib_b^* v}+\not{\NO{\dib_b v}}+\not{\NO{v}_{-1}}
\\&=\NO{u}_0.
\end{split}
\end{equation*}
We use the microlocal decompositin of $u$ and $v$ and notice that $u^+=\dib_b^*( v^+)+[\dib_b^*,\Psi^+]v$ and that $[\dib_b^*,\Psi^+]v=\tilde\Psi^0 v$. Hence
\begin{equation*}
\begin{split}
\NO{u^+}&\leq (u^+,\dib_b^*( v^+))+(u^+,\tilde\Psi^0 v)
\\
&\leq (\dib_b (u^+),v^+)+(\tilde\Psi^0u,v)
\\
&\simleq \no{\dib_b(u^+)}\epsilon(\no{\dib_b(v^+)}+\no{(\dib_b^*v)^+}+\no{\tilde \Psi^0v})+\no{\tilde\Psi^0u}\no{u}
\\
&\leq \no{\dib_b(u^+)}\epsilon\no{u^+}+\epsilon\no{\dib_b(u^+)}\no{\tilde\Psi^0v}+\no{u}\no{\tilde\Psi^0u}
\\
&\leq \epsilon\NO{\dib_b(u^+)}+\epsilon\no{u^+}^2+s.c.\no{u}^2+l.c.\no{\tilde\Psi^0u}^2,
\end{split}
\end{equation*}
where s.c. and l.c. denote a small and large constant respectively. If we take summation over $+,-,0$, we end up with
$$
\no{u}\leq\epsilon\no{\dib_b u}+l.c.\no{\tilde\Psi^0u}.
$$
When $n=2$ and only estimates for $v^+$ and not for the full $v$ in degree $1$  are  provided by Theorem~\ref{t1.1}, we use the extra assumption that $\dib_b$ has closed range. Thus $\dib_b^*$ has also closed range, we can write $u\in\H^{0\,\perp}$ as $u=\dib_b^* v$ for $v$ satisfying $\dib_bv=0$ and $\no{v}\simleq \no{u}$ and the proof goes through in the same way as above by using \eqref{2.2} for $u^+$.

\hskip15cm$\Box$

\bibliographystyle{alphanum}

\end{document}